\documentclass{article}
\usepackage[utf8]{inputenc}
\usepackage{amsmath}
\usepackage{amssymb}
\usepackage{breqn}

\title{A General Formula for Asymptotes of Rational Polynomial Functions}
\author{\quad\quad\quad Lam Mason, \quad\quad\quad Asterios Skodras \\   mason.Lam@ccoex.com, \quad askodras@sch.gr}

\date{April 2021}

\begin{document}

\maketitle

{\abstract}
We propose a formula for finding the horizontal, oblique or curvilinear asymptote of any rational polynomial function of any positive degree, as a sum of matrix determinants formed directly from the coefficients of the terms in the given polynomial. This formula provides a new means of computing asymptotes in addition to the standard methods of Euclidean division and the evaluation of limits. \\
\\
\providecommand{\keywords}[1]
{
  \small
  \textbf{\text{Keywords:}} #1
}
\( \\ \\ \)
\section {Introduction}
It is well known, and explained in many textbooks [1, 2, 3, 4], that the \textit{horizontal}, \textit{oblique} or \textit{curvilinear} asymptote of a function \(f:D \to f(D), D\subseteq \mathrm R\) is a straight line or a curve whose distance from the graph of \(f(x)\) approaches zero as \(x\) approaches infinity. In the case of a rational polynomial function \(\frac{a(x)}{b(x)}\), where deg(\(b(x))\leq \deg (a(x))\), the expression for the asymptote is a polynomial with degree deg\(\left(a(x)\right)\)-deg\(\left(b(x)\right)\). Currently one of the most common methods of finding the asymptote of a rational polynomial function is by \textit{Euclidean polynomial division}, where the numerator of the polynomial is divided by denominator and the fractional terms are discarded. \\
\\
\textbf{Example 1.1} Consider the function
\begin{dmath}
    f(x)=\frac{8x^3+7}{x-4}
\end{dmath}
By Euclidean division,
\begin{dmath}
    f(x)=8x^2+32x+128+\frac{519}{x-4}
\end{dmath}
Hence, the asymptote of this function is
\begin{dmath}
    g(x)=8x^2+32x+128
\end{dmath}
Another popular method, specifically used for finding \textit{oblique} asymptotes, is by evaluating limits. The oblique asymptote of a function \(f(x)\) takes the form
\begin{dmath}
    y=ax+b
\end{dmath}
where
\begin{dmath}
    a=\lim_{x\to\infty}\frac{f(x)}{x}
\end{dmath}
and
\begin{dmath}
    b=\lim_{x\to\infty}\big(f(x)-ax\big)
\end{dmath}
\textbf{Example 1.2} For the function
\begin{dmath}
    f(x)=\frac{5x^3+13x^2+3x+9}{4x^2+5x+7}
\end{dmath}
we have
\begin{dmath}
    a=\lim_{x\to\infty}\frac{5x^3+13x^2+3x+9}{4x^3+5x^2+7x}
\end{dmath}
\begin{dmath}
    a=\frac{5}{4}
\end{dmath}
and
\begin{dmath}
    b=\lim_{x\to\infty}\frac{5x^3+13x^2+3x+9}{4x^2+5x+7}-\frac{5}{4}x \\
\end{dmath}
\begin{dmath}
    b=\frac{27}{16}
\end{dmath}
Thus, the corresonding oblique asymptote is given by
\begin{dmath}
    g(x)=\frac{5}{4}x+\frac{27}{16}
\end{dmath}
In this paper, we show that the asymptote can alternatively be expressed as the following sum:
\begin{dmath}
    g(x)=\frac{a_n}{b_{n-k}}x^k-\frac{
    \begin{vmatrix}
        a_n & a_{n-1} \\
        b_{n-k} & b_{n-k-1}
    \end{vmatrix}}{b_{n-k}^2}x^{k-1}+\frac{
    \begin{vmatrix}
        a_n & a_{n-1} & a_{n-2} \\
        b_{n-k} & b_{n-k-1} & b_{n-k-2} \\
        0 & b_{n-k} & b_{n-k-1}
    \end{vmatrix}}{b_{n-k}^3}x^{k-2}+\dots+(-1)^{k-1}\frac{
    \begin{vmatrix}
        a_n & a_{n-1} & \dots & a_{n-(k-1)} \\
        b_{n-k} & b_{n-k-1} & \dots & b_{n-k-(k-1)} \\
        0 & b_{n-k} & \dots & b_{n-k-(k-1)+1} \\
        \vdots & \vdots & \ddots & \vdots \\
        0 & 0 & \dots & b_{n-k-1}
    \end{vmatrix}}{b_{n-k}^k}x+(-1)^k\frac{
    \begin{vmatrix}
        a_n & a_{n-1} & \dots & a_{n-k} \\
        b_{n-k} & b_{n-k-1} & \dots & b_{n-k-k} \\
        0 & b_{n-k} & \dots & b_{n-k-k+1} \\
        \vdots & \vdots & \ddots & \vdots \\
        0 & 0 & \dots & b_{n-k-1}
    \end{vmatrix}}{b_{n-k}^{k+1}}
\end{dmath}
where \(b_i=0\) for \(i<0\), for any rational polynomial function
\begin{dmath}
    f(x)=\frac{a_nx^n+a_{n-1}x^{n-1}+\dots+a_1x+a_0}{b_{n-k}x^{n-k}+b_{n-k-1}x^{n-k-1}+\dots+b_1x+b_0}
\end{dmath}
where \(n>k\) and \(n, k\in\mathbb{Z}\).
\( \\ \\ \)
\section{Proof by Induction}
For a rational polynomial function
\begin{equation}
   f(x)= \displaystyle \frac{\displaystyle\displaystyle\sum_{r=0}^na_rx^r}{\displaystyle\sum_{i=0}^{n-k} b_ix^i}
\end{equation}
the equation for its corresponding asymptote would be in the form of
\begin{equation}
    g(x)=\displaystyle\sum_{\alpha=0}^k\theta_\alpha x^{k-\alpha}
\end{equation}
for some coefficients \(\theta_\alpha\). \\
We first derive the expression for \(\theta_0\), the coefficient of \(x^k\). By the definition of an asymptote,
\begin{equation}
    \lim_{x\to\infty}\left(f(x)-g(x)\right)=0
\end{equation}
Plugging in our expression for \(g(x)\) gives
\begin{equation}
    \lim_{x\to\infty}\left(f(x)-\displaystyle\sum_{\alpha=0}^k\theta_\alpha x^{k-\alpha}\right)=0
\end{equation}
\begin{equation}
    \lim_{x\to\infty}\left(f(x)-\displaystyle\sum_{\alpha=1}^k\theta_\alpha x^{k-\alpha}-\theta_0x^k\right)=0
\end{equation}
Dividing both sides by \(x^k\),
\begin{equation}
     \lim_{x\to\infty}\left(\frac{f(x)}{x^k}-\displaystyle\sum_{\alpha=1}^k\theta_\alpha x^{-\alpha}-\theta_0\right)=0
\end{equation}
Since
\begin{equation}
    \displaystyle\sum_{\alpha=1}^k\theta_\alpha x^{-\alpha}
\end{equation}
tends to 0 as \(x\to\infty\), we have
\begin{equation}
    \lim_{x\to\infty}\left(\frac{f(x)}{x^k}-\theta_0\right)=0
\end{equation}
\begin{equation}
    \lim_{x\to\infty}\left(\frac{f(x)}{x^k}\right)-\theta_0=0
\end{equation}
\begin{equation}
    \lim_{x\to\infty}\left(\frac{f(x)}{x^k}\right)=\theta_0
\end{equation}
Expanding \(f(x)\), we have
\begin{equation}
    \theta_0=\lim_{x\to\infty}\frac{f(x)}{x^k}=\lim_{x\to\infty}\frac{\displaystyle\sum_{r=0}^{n}a_rx^r}{x^k\displaystyle\sum_{i=0}^{n-k}b_ix^i}=\lim_{x \to \infty}\frac{\displaystyle\sum_{r=0}^na_rx^r}{\displaystyle\sum_{i=k}^nb_{i-k}x^i}
\end{equation}
\begin{equation}
    \theta_0=\lim_{x\to\infty}\frac{a_nx^n}{b_{n-k}x^n}=\lim_{x\to\infty}\frac{a_n}{b_{n-k}}=\frac{a_n}{b_{n-k}}
\end{equation}
Now, we prove that
\begin{equation}
    \theta_\alpha=\left(-1\right)^{\alpha}\frac{
    \begin{vmatrix}
        a_n & a_{n-1} & \dots & a_{n-\alpha} \\
        b_{n-k} & b_{n-k-1} & \dots & b_{n-k-\alpha} \\
        0 & b_{n-k} & \dots & b_{n-k-\alpha+1} \\
        \vdots & \vdots & \ddots & \vdots \\
        0 & 0 & \dots & b_{n-k-1}
    \end{vmatrix}}{b_{n-k}^{\alpha+1}}
\end{equation}
where \(\alpha\leq k\in\mathbb{Z}^+\), \(b_i=0\) for \(i<0\), and we do so by strong induction.\\
\\
Let \(P(\alpha)\) be the proposition stated above. Then, to find \(\theta_\alpha\) for \(\alpha=1\), we begin with
\begin{equation}
    \lim_{x\to\infty}\left(f(x)-\displaystyle\sum_{\alpha=0}^k\theta_\alpha x^{k-\alpha}\right)=0
\end{equation}
\begin{equation}
    \lim_{x\to\infty}\left(f(x)-\displaystyle\sum_{\alpha=2}^k\theta_\alpha x^{k-\alpha}-\theta_1x^{k-1}-\theta_0x^k\right)=0
\end{equation}
Dividing both sides by \(x^{k-1}\),
\begin{equation}
    \lim_{x\to\infty}\left(\frac{f(x)}{x^{k-1}}-\displaystyle\sum_{\alpha=2}^{k}\theta_\alpha x^{-\alpha+1}-\theta_1-\theta_0x\right)=0
\end{equation}
\begin{equation}
    \lim_{x\to\infty}\left(\frac{f(x)}{x^{k-1}}-\theta_1-\theta_0x\right)=0
\end{equation}
\begin{equation}
    \lim_{x\to\infty}\left(\frac{f(x)}{x^{k-1}}-\theta_0x\right)=0
\end{equation}
\begin{equation}
    \lim_{x\to\infty}\left(\frac{f(x)}{x^{k-1}}-\theta_0x\right)-\theta_1=0
\end{equation}
\begin{equation}
    \lim_{x\to\infty}\left(\frac{f(x)}{x^{k-1}}-\theta_0x\right)=\theta_1
\end{equation}
Expanding \(f(x)\) and referring back to our definition of \(\theta_0\) yields
\begin{equation}
    \theta_1=\lim_{x\to\infty}\left(\frac{f(x)}{x^{k-1}}-\theta_0x\right)=\lim_{x\to\infty}\left(\frac{\displaystyle\sum_{r=0}^{n}a_rx^r}{x^{k-1}\displaystyle\sum_{i=0}^{n-k}b_ix^i}-\frac{a_n}{b_{n-k}}x\right)
\end{equation}
\begin{equation}
    \theta_1=\lim_{x\to\infty}\left(\frac{\displaystyle\sum_{r=0}^na_rb_{n-k}x^r-\displaystyle\sum_{i=k}^na_nb_{i-k}x^i}{\displaystyle\sum_{i=k-1}^{n-1}b_{i-k+1}b_{n-k}x^i}\right)
\end{equation}
\begin{equation}
    \theta_1=\lim_{x\to\infty}\frac{\left(a_rb_{n-k}x^r-a_nb_{n-k-1}\right)x^{n-1}}{b_{n-k}{b_{n-k}x^{n-1}}}=\lim_{x\to\infty}\frac{a_{n-1}b_{n-k}-a_nb_{n-k-1}}{b_{n-k}^2}
\end{equation}
\begin{equation}
    \theta_1=\frac{a_{n-1}b_{n-k}-a_nb_{n-k-1}}{b_{n-k}^2}=-\frac{a_nb_{n-k-1}-a_{n-1}b_{n-k}}{b_{n-k}^2}
\end{equation}
\begin{equation}
    \theta_1=-\frac{
    \begin{vmatrix}
        a_n & a_{n-1} \\
        b_{n-k} & b_{n-k-1} \\
    \end{vmatrix}}{b_{n-k}^2}
\end{equation}
Thus, \(P(\alpha)\) is true.
Now, suppose \(P(\alpha)\) is true for all \(\alpha\leq\varphi-1\leq k-1\in\mathbb{Z}^+\). We will then show that \(P(\varphi)\) is true, i.e.
\begin{equation}
    \theta_\varphi=\left(-1\right)^\varphi\frac{
    \begin{vmatrix}
        a_n & a_{n-1} & \dots & a_{n-\varphi} \\
        b_{n-k} & b_{n-k-1} & \dots & b_{n-k-\varphi} \\
        0 & b_{n-k} & \dots & b_{n-k-\varphi+1} \\
        \vdots & \vdots & \ddots & \vdots \\
        0 & 0 & \dots & b_{n-k-1}
    \end{vmatrix}}{b_{n-k}^{\varphi+1}}
\end{equation}
for \(\varphi-1<\varphi\leq k\in\mathbb{Z}^+\). \\
First we have
\begin{equation}
    \theta_\varphi=\lim_{x\to\infty}\left(\frac{f(x)}{x^{k-\varphi}}-\displaystyle\sum_{\alpha=0}^{\varphi-1}\theta_\alpha x^{\varphi-\alpha}\right)
\end{equation}
\begin{dmath}
    \theta_\varphi=\lim_{x\to\infty}\left(\frac{f(x)}{x^{k-\varphi}}-\frac{a_n}{b_{n-k}}x^\varphi-(-1)\frac{
    \begin{vmatrix}
        a_n & a_{n-1} \\
        b_{n-k} & b_{n-k-1} \\
    \end{vmatrix}}{b_{n-k}^2}x^{\varphi-1}-\dots-(-1)^{\varphi-1}\frac{
    \begin{vmatrix}
        a_n & \dots & a_{n-(\varphi-1)} \\
        b_{n-k} & \dots & b_{n-k-(\varphi-1)} \\
        \vdots & \ddots & \vdots \\
        0 & \dots & b_{n-k-1}
    \end{vmatrix}}{b_{n-k}^\varphi}x\right)
\end{dmath}
\begin{dmath}
    \theta_\varphi=\lim_{x\to\infty}\left(\frac{\displaystyle\sum_{r=0}^na_rx^r}{x^{k-\varphi}\displaystyle\sum_{i=0}^{n-k}b_ix^i}-\frac{a_n}{b_{n-k}}x^\varphi-\dots-\frac{(-1)^{\varphi-1}
    \begin{vmatrix}
        a_n & \dots & a_{n-{\varphi-1}} \\
        b_{n-k} & \dots & b_{n-{\varphi-1}} \\
        \vdots & \ddots & \vdots \\
        0 & \dots & b_{n-k-1}
    \end{vmatrix}}{b_{n-k}^\varphi}x\right)
\end{dmath}
Let
\begin{equation}
    \theta_\alpha^{'}=(-1)^\alpha
    \begin{vmatrix}
        a_n & \dots & a_{n-\alpha} \\
        b_{n-k} & \dots & b_{n-k-\alpha} \\
        \vdots & \ddots & \vdots \\
        0 & \dots & b_{n-k-1}
    \end{vmatrix}
\end{equation}
Then, combining the fractions gives
\begin{equation}
    \theta_\varphi=\lim_{x\to\infty}\left(\frac{\displaystyle\sum_{r=0}^na_rx^r}{x^{k-\varphi}\displaystyle\sum_{i=0}^{n-k}b_ix^i}-\frac{\displaystyle\sum_{\alpha=0}^{\varphi-1}b_{n-k}^{\varphi-1-\alpha}\theta_\alpha^{'}x^{\varphi-\alpha}}{b_{n-k}^\varphi}\right)
\end{equation}
\begin{equation}
    \theta_\varphi=\lim_{x\to\infty}\frac{\displaystyle\sum_{r=0}^na_rb_{n-k}^\varphi x^r-\left(\displaystyle\sum_{i=0}^{n-k}b_ix^i\right)\left(\displaystyle\sum_{\alpha=0}^{\varphi-1}b_{n-k}^{\varphi-1-\alpha}\theta_\alpha^{'}x^{\varphi-\alpha}\right)}{\displaystyle\sum_{i=k}^{n-\varphi}b_{i-k+\varphi}b_{n-k}^\varphi x^i}
\end{equation}
We show that all the terms in the numerator where the powers of \(x\) are greater than \(n-\varphi\) vanish using the following result: \\
\\
\textbf{Lemma 2.1} In the expression
\begin{dmath}
    \sum_{r=0}^na_rb_{n-k}^\varphi x^r-\left(\displaystyle\sum_{i=0}^{n-k}b_ix^i\right)\left(\displaystyle\sum_{\alpha=0}^{\varphi-1}b_{n-k}^{\varphi-1-\alpha}\theta_\alpha^{'}x^{\varphi-\alpha}\right)
\end{dmath}
the terms where the powers of \(x\) are greater than \(n-\varphi\) all sum to zero. \\
\\
By this result, we are left with
\begin{dmath}
    \theta_\varphi=\lim_{x\to\infty}\frac{1}{b_{n-k}b_{n-k}^\varphi x^{n-\varphi}}\left(a_{n-\varphi}b_{n-k}^\varphi-\left(b_{n-k-\varphi}b_{n-k}^{\varphi-1}\theta_0^{'}+b_{n-k-(\varphi-1)}b_{n-k}^{\varphi-2}\theta_1^{'}+\dots+b_{n-k-2}b_{n-k}\theta_{\varphi-2}^{'}+b_{n-k-1}\theta_{\varphi-1}^{'}\right)\right)x^{n-\varphi}
\end{dmath}
\begin{dmath}
    \theta_\varphi=\lim_{x\to\infty}\frac{1}{b_{n-k}^{\varphi+1}}\left(a_{n-\varphi}b_{n-k}^\varphi-\left(b_{n-k-\varphi}b_{n-k}^{\varphi-1}\theta_0^{'}+b_{n-k-(\varphi-1)}b_{n-k}^{\varphi-2}\theta_1^{'}+\dots+b_{n-k-2}b_{n-k}\theta_{\varphi-2}^{'}+b_{n-k-1}\theta_{\varphi-1}^{'}\right)\right)
\end{dmath}
To simplify this expression for \(\theta_\varphi\), we make use of the following result: \\
\\
\textbf{Lemma 2.2}
\begin{dmath}
    \theta_j^{'}=a_{n-j}b_{n-k}^j-\left(b_{n-k-j}b_{n-k}^{j-1}\theta_0^{'}+b_{n-k-(j-1)}b_{n-k}^{j-2}\theta_1^{'}+\dots+b_{n-k-2}b_{n-k}\theta_{j-2}^{'}+b_{n-k-1}\theta_{j-1}^{'}\right)
\end{dmath}
which gives us
\begin{dmath}
    \theta_\varphi^{'}=a_{n-\varphi}b_{n-k}^\varphi-(b_{n-k-\varphi}b_{n-k}^{\varphi-1}\theta_0^{'}+b_{n-k-(\varphi-1)}b_{n-k}^{\varphi-2}\theta_1^{'}+\dots+b_{n-k-2}b_{n-k}\theta_{\varphi-2}^{'}+b_{n-k-1}\theta_{\varphi-1}^{'})
\end{dmath}
Hence, we have
\begin{dmath}
    \theta_\varphi^{'}=\frac{\theta_\varphi^{'}}{b_{n-k}^{\varphi+1}}
\end{dmath}
\begin{dmath}
    \theta_\varphi^{'}=(-1)^\varphi\frac{
    \begin{vmatrix}
        a_n & a_{n-1} & \dots & a_{n-(\varphi-2)} & a_{n-(\varphi-1)} & a_{n-\varphi} \\
        b_{n-k} & b_{n-k-1} & \dots & b_{n-k-(\varphi-2)} & b_{n-k-(\varphi-1)} & b_{n-k-\varphi} \\
        \vdots & \vdots & \ddots & \vdots & \vdots & \vdots \\
        0 & 0 & \dots & b_{n-k} & b_{n-k-1} & b_{n-k-2} \\
        0 & 0 & \dots & 0 & b_{n-k} & b_{n-k-1}
    \end{vmatrix}}{b_{n-k}^{\varphi+1}}
\end{dmath}
Hence, \(P(\varphi)\) is true. \\
This means \(P(\varphi)\) is true for all \(\alpha\leq k\in\mathbb{Z}^+\), thus completing the proof.
\( \\ \\ \)
\section{Proof of Lemmas}
In this section we prove the two lemmas mentioned previously.
\subsection{Proof of Lemma 2.1}
The terms in the expression
\begin{dmath}
    \left(\displaystyle\sum_{i=0}^{n-k}b_ix^i\right)\left(\displaystyle\sum_{a=0}^{\varphi-1}b_{n-k}^{\varphi-1-\alpha}\theta_\alpha^{'}x^{k-\alpha}\right)
\end{dmath}
where the powers of \(x\) are greater than \(n-\varphi\) are
\begin{dmath}
    \left(b_{n-k-0}x^{n-k-0}\right)\left(b_{n-k}^{\varphi-1-z}\theta_z^{'}x^{k-z}\right)+\left(b_{n-k-1}x^{n-k-1}\right)\left(b_{n-k}^{\varphi-1-(z-1)}\theta_{z-1}^{'}x^{k-(z-1)}\right)+\left(b_{n-k-2}x^{n-k-2}\right)\left(b_{n-k}^{\varphi-1-(z-2)}\theta_{z-2}^{'}x^{k-(z-2)}\right)+\dots+\left(b_{n-k-z}x^{n-k-z}\right)\left(b_{n-k}^{\varphi-1-0}\theta_0^{'}x^{k-0}\right)
\end{dmath}
for integers \(0\leq z\leq \varphi-1\). Thus, in order to prove that the terms in the expression
\begin{dmath}
    \sum_{r=0}^na_rb_{n-k}^\varphi x^r-\left(\sum_{i=0}^{n-k}b_ix^i\right)\left(\sum_{\alpha=0}^{\varphi-1}b_{n-k}^{\varphi-1-\alpha}\theta_\alpha^{'}x^{k-\alpha}\right)
\end{dmath}
all sum to zero for powers of \(x\) greater than \(n-\varphi\), we need to show that
\begin{dmath}
    a_{n-z}b_{n-k}^\varphi x_{n-z}=\left(b_{n-k-0}x^{n-k-0}\right)\left(b_{n-k}^{\varphi-1-z}\theta_z^{'}x^{k-z}\right)+\left(b_{n-k-1}x^{n-k-1}\right)\left(b_{n-k}^{\varphi-1-(z-1)}\theta_{z-1}^{'}x^{k-(z-1)}\right)+\left(b_{n-k-2}x^{n-k-2}\right)\left(b_{n-k}^{\varphi-1-(z-2)}\theta_{z-2}^{'}x^{k-(z-2)}\right)+\dots+\left(b_{n-k-z}x^{n-k-z}\right)\left(b_{n-k}^{\varphi-1-0}\theta_0^{'}x^{k-0}\right)
\end{dmath}
i.e.
\begin{dmath}
    a_{n-z}b_{n-k}^\varphi=\left(b_{n-k-0}x^{n-k-0}\right)\left(b_{n-k}^{\varphi-1-z}\theta_z^{'}\right)+\left(b_{n-k-1}\right)\left(b_{n-k}^{\varphi-1-(z-1)}\theta_{z-1}^{'}\right)+\left(b_{n-k-2}\right)\left(b_{n-k}^{\varphi-1-(z-2)}\theta_{z-2}^{'}\right)+\dots+\left(b_{n-k-z}\right)\left(b_{n-k}^{\varphi-1-0}\theta_0^{'}\right)
\end{dmath}
To do this, we refer to Lemma 2.2 for our expansion of \(\theta_z^{'}\), which gives us
\begin{dmath}
    \theta_z^{'}=a_{n-z}b_{n-k}^z-\left(b_{n-k-z}b_{n-k}^{z-1}\theta_0^{'}+b_{n-k-(z-1)}b_{n-k}^{z-2}\theta_1^{'}+\dots+b_{n-k-2}b_{n-k}\theta_{z-2}^{'}+b_{n-k-1}\theta_{z-1}^{'}\right)
\end{dmath}
Plugging this expansion onto the RHS gives
\begin{dmath}
    \text{RHS}=\left(b_{n-k-0}\right)\left\{b_{n-k}^{\varphi-1-z}\left[a_{n-z}b_{n-k}^z-\left(b_{n-k-z}b_{n-k}^{z-1}\theta_0^{'}+b_{n-k-(z-1)}b_{n-k}^{z-2}\theta_1^{'}+\dots+b_{n-k-2}b_{n-k}\theta_{z-2}^{'}+b_{n-k-1}\theta_{z-1}^{'}\right)\right]\right\}+\left(b_{n-k-1}\right)\left(b_{n-k}^{\varphi-1-(z-1)}\theta_{z-1}^{'}\right)+\left(b_{n-k-2}\right)\left(b_{n-k}^{\varphi-1-(z-2)}\theta_{z-2}^{'}\right)+\dots+\left(b_{n-k-z}\right)\left(b_{n-k}^{\varphi-1-0}\theta_0^{'}\right)
\end{dmath}
\begin{dmath}
    \text{RHS}=a_{n-z}b_{n-k}^\varphi-\left(b_{n-k-z}b_{n-k}^{\varphi-1}\theta_0^{'}+b_{n-k-(z-1)}b_{n-k}^{\varphi-2}\theta_1^{'}+\dots+b_{n-k-2}b_{n-k}^{\varphi-(z-1)}\theta_{z-2}^{'}+b_{n-k-1}b_{n-k}^{\varphi-z}\theta_{z-1}^{'}\right)+\left(b_{n-k-1}\right)\left(b_{n-k}^{\varphi-1-(z-1)}\theta_{z-1}^{'}\right)+\left(b_{n-k-2}\right)\left(b_{n-k}^{\varphi-1-(z-2)}\theta_{z-2}^{'}\right)+\dots+\left(b_{n-k-z}\right)\left(b_{n-k}^{\varphi-1-0}\theta_0^{'}\right)
\end{dmath}
\begin{dmath}
    \text{RHS}=a_{n-z}b_{n-k}^\varphi-\left(b_{n-k-z}b_{n-k}^{\varphi-1}\theta_0^{'}+b_{n-k-(z-1)}b_{n-k}^{\varphi-2}\theta_1^{'}+\dots+b_{n-k-2}b_{n-k}^{\varphi-(z-1)}\theta_{z-2}^{'}+b_{n-k-1}b_{n-k}^{\varphi-z}\theta_{z-1}^{'}\right)+\left(b_{n-k-1}\right)\left(b_{n-k}^{\varphi-z}\theta_{z-1}^{'}\right)+\left(b_{n-k-2}\right)\left(b_{n-k}^{\varphi-(z-1)}\theta_{z-2}^{'}\right)+\dots+\left(b_{n-k-z}\right)\left(b_{n-k}^{\varphi-1}\theta_0^{'}\right)
\end{dmath}
\begin{dmath}
    \text{RHS}=a_{n-z}b_{n-k}^{\varphi}
\end{dmath}
as desired.
\subsection{Proof of Lemma 2.2}
We begin with
\begin{equation}
    \theta_j^{'}=(-1)^j
    \begin{vmatrix}
        a_n & a_{n-1} & \dots & a_{n-(j-2)} & a_{n-(j-1)} & a_{n-j} \\
        b_{n-k} & b_{n-k-1} & \dots & b_{n-k-(j-2)} & b_{n-k-(j-1)} & b_{n-k-j} \\
        \vdots & \vdots & \ddots & \vdots & \vdots & \vdots \\
        0 & 0 & \dots & b_{n-k} & b_{n-k-1} & b_{n-k-2} \\
        0 & 0 & \dots & 0 & b_{n-k} & b_{n-k-1}
    \end{vmatrix}
\end{equation}
Dividing both sides by \((-1)^j\),
\begin{equation}
    \frac{\theta_j^{'}}{(-1)^j}=
    \begin{vmatrix}
        a_n & a_{n-1} & \dots & a_{n-(j-2)} & a_{n-(j-1)} & a_{n-j} \\
        b_{n-k} & b_{n-k-1} & \dots & b_{n-k-(j-2)} & b_{n-k-(j-1)} & b_{n-k-j} \\
        \vdots & \vdots & \ddots & \vdots & \vdots & \vdots \\
        0 & 0 & \dots & b_{n-k} & b_{n-k-1} & b_{n-k-2} \\
        0 & 0 & \dots & 0 & b_{n-k} & b_{n-k-1}
    \end{vmatrix}
\end{equation}
Expanding the determinant by taking the bottom row elements gives
\begin{dmath}
    \frac{\theta_j^{'}}{(-1)^j}=-b_{n-k}
    \begin{vmatrix}
        a_n & a_{n-1} & \dots & a_{n-(j-2)} & a_{n-j} \\
        b_{n-k} & b_{n-k-1} & \dots & b_{n-k-(j-2)} & b_{n-k-j} \\
        \vdots & \vdots & \ddots & \vdots & \vdots \\
        0 & 0 & \dots & b_{n-k-1} & b_{n-k-3} \\
        0 & 0 & \dots & b_{n-k} & b_{n-k-2}
    \end{vmatrix}+b_{n-k-1}
    \begin{vmatrix}
        a_n & a_{n-1} & \dots & a_{n-(j-2)} & a_{n-(j-1)} \\
        b_{n-k} & b_{n-k-1} & \dots & b_{n-k-(j-2)} & b_{n-k-(j-1)} \\
        \vdots & \vdots & \ddots & \vdots & \vdots \\
        0 & 0 & \dots & b_{n-k-1} & b_{n-k-2} \\
        0 & 0 & \dots & b_{n-k} & b_{n-k-1}
    \end{vmatrix}
\end{dmath}
\begin{dmath}
    \frac{\theta_j^{'}}{(-1)^j}=-b_{n-k}
    \begin{vmatrix}
        a_n & a_{n-1} & \dots & a_{n-(j-2)} & a_{n-j} \\
        b_{n-k} & b_{n-k-1} & \dots & b_{n-k-(j-2)} & b_{n-k-j} \\
        \vdots & \vdots & \ddots & \vdots & \vdots \\
        0 & 0 & \dots & b_{n-k-1} & b_{n-k-3} \\
        0 & 0 & \dots & b_{n-k} & b_{n-k-2}
    \end{vmatrix}+b_{n-k-1}\frac{\theta_{j-1}^{'}}{(-1)^{j-1}}
\end{dmath}
Expanding the determinant again, we get
\begin{dmath}
     \frac{\theta_j^{'}}{(-1)^j}=-b_{n-k}\left(-b_{n-k}
     \begin{vmatrix}
        a_n & a_{n-1} & \dots & a_{n-(j-3)} & a_{n-j} \\
        b_{n-k} & b_{n-k-1} & \dots & b_{n-k-(j-3)} & b_{n-k-j} \\
        \vdots & \vdots & \ddots & \vdots & \vdots \\
        0 & 0 & \dots & b_{n-k-1} & b_{n-k-4} \\
        0 & 0 & \dots & b_{n-k} & b_{n-k-3}
     \end{vmatrix}+b_{n-k-2}
     \begin{vmatrix}
        a_n & a_{n-1} & \dots & a_{n-(j-3)} & a_{n-(j-2)} \\
        b_{n-k} & b_{n-k-1} & \dots & b_{n-k-(j-3)} & b_{n-k-(j-2)} \\
        \vdots & \vdots & \ddots & \vdots & \vdots \\
        0 & 0 & \dots & b_{n-k-1} & b_{n-k-2} \\
        0 & 0 & \dots & b_{n-k} & b_{n-k-3}
     \end{vmatrix}\right)+b_{n-k-1}\frac{\theta_{j-1}^{'}}{(-1)^{j-1}}
\end{dmath}
\begin{dmath}
    \frac{\theta_j^{'}}{(-1)^j}=-b_{n-k}\left(-b_{n-k}
    \begin{vmatrix}
        a_n & a_{n-1} & \dots & a_{n-(j-3)} & a_{n-j} \\
        b_{n-k} & b_{n-k-1} & \dots & b_{n-k-(j-3)} & b_{n-k-j} \\
        \vdots & \vdots & \ddots & \vdots & \vdots \\
        0 & 0 & \dots & b_{n-k-1} & b_{n-k-4} \\
        0 & 0 & \dots & b_{n-k} & b_{n-k-3}
     \end{vmatrix}+b_{n-k-2}\frac{\theta_{j-2}^{'}}{(-1)^{j-2}}\right)+b_{n-k-1}\frac{\theta_{j-1}^{'}}{(-1)^{j-1}}
\end{dmath}
\begin{dmath}
    \frac{\theta_j^{'}}{(-1)^j}=b_{n-k}^2
    \begin{vmatrix}
        a_n & a_{n-1} & \dots & a_{n-(j-3)} & a_{n-j} \\
        b_{n-k} & b_{n-k-1} & \dots & b_{n-k-(j-3)} & b_{n-k-j} \\
        \vdots & \vdots & \ddots & \vdots & \vdots \\
        0 & 0 & \dots & b_{n-k-1} & b_{n-k-4} \\
        0 & 0 & \dots & b_{n-k} & b_{n-k-3}
    \end{vmatrix}-b_{n-k}b_{n-k-2}\frac{\theta_{j-2}^{'}}{(-1)^{j-2}}+b_{n-k-1}\frac{\theta_{j-1}^{'}}{(-1)^{j-1}}
\end{dmath}
Expanding once again gives
\begin{dmath}
    \frac{\theta_j^{'}}{(-1)^j}=b_{n-k}^2\left(-b_{n-k}
    \begin{vmatrix}
        a_n & a_{n-1} & \dots & a_{n-(j-4)} & a_{n-j} \\
        b_{n-k} & b_{n-k-1} & \dots & b_{n-k-(j-4)} & b_{n-k-j} \\
        \vdots & \vdots & \ddots & \vdots & \vdots \\
        0 & 0 & \dots & b_{n-k-1} & b_{n-k-5} \\
        0 & 0 & \dots & b_{n-k} & b_{n-k-4}
    \end{vmatrix}+b_{n-k-3}
    \begin{vmatrix}
        a_n & a_{n-1} & \dots & a_{n-(j-4)} & a_{n-(j-3)} \\
        b_{n-k} & b_{n-k-1} & \dots & b_{n-k-(j-4)} & b_{n-k-(j-3)} \\
        \vdots & \vdots & \ddots & \vdots & \vdots \\
        0 & 0 & \dots & b_{n-k-1} & b_{n-k-2} \\
        0 & 0 & \dots & b_{n-k} & b_{n-k-1}
    \end{vmatrix}\right)-b_{n-k}b_{n-k-2}\frac{\theta_{j-2}^{'}}{(-1)^{j-2}}+b_{n-k-1}\frac{\theta_{j-1}^{'}}{(-1)^{j-1}}
\end{dmath}
\begin{dmath}
    \frac{\theta_j^{'}}{(-1)^j}=b_{n-k}^2\left(-b_{n-k}
    \begin{vmatrix}
        a_n & a_{n-1} & \dots & a_{n-(j-4)} & a_{n-j} \\
        b_{n-k} & b_{n-k-1} & \dots & b_{n-k-(j-4)} & b_{n-k-j} \\
        \vdots & \vdots & \ddots & \vdots & \vdots \\
        0 & 0 & \dots & b_{n-k-1} & b_{n-k-5} \\
        0 & 0 & \dots & b_{n-k} & b_{n-k-4}
    \end{vmatrix}+b_{n-k-2}\frac{\theta_{j-3}^{'}}{(-1)^{j-3}}\right)-b_{n-k}b_{n-k-2}\frac{\theta_{j-2}^{'}}{(-1)^{j-2}}+b_{n-k-1}\frac{\theta_{j-1}^{'}}{(-1)^{j-1}}
\end{dmath}
\begin{dmath}
    \frac{\theta_j^{'}}{(-1)^j}=-b_{n-k}^3
    \begin{vmatrix}
        a_n & a_{n-1} & \dots & a_{n-(j-4)} & a_{n-j} \\
        b_{n-k} & b_{n-k-1} & \dots & b_{n-k-(j-4)} & b_{n-k-j} \\
        \vdots & \vdots & \ddots & \vdots & \vdots \\
        0 & 0 & \dots & b_{n-k-1} & b_{n-k-5} \\
        0 & 0 & \dots & b_{n-k} & b_{n-k-4}
    \end{vmatrix}+b_{n-k}^2{n-k-2}\frac{\theta_{j-3}^{'}}{(-1)^{j-3}}-b_{n-k}b_{n-k-2}\frac{\theta_{j-2}^{'}}{(-1)^{j-2}}+b_{n-k-1}\frac{\theta_{\varphi-1}^{'}}{(-1)^{j-1}}
\end{dmath}
Repeating this process yields
\begin{dmath}
    \frac{\theta_j^{'}}{(-1)^j}=(-1)^{j-1}b_{n-k}^{j-1}
    \begin{vmatrix}
        a_n & a_{n-j} \\
        b_{n-k} & b_{n-k-j}
    \end{vmatrix}+(-1)^{j-2}b_{n-k}^{j-2}b_{n-k-(j-1)}\frac{\theta_1^{'}}{(-1)^1}+(-1)^{j-3}b_{n-k}^{j-3}b_{n-k-(j-2)}\frac{\theta_2^{'}}{(-1)^2}+\dots+(-1)^2b_{n-k}^2b_{n-k-3}\frac{\theta_{j-3}^{'}}{(-1)^{j-3}}+(-1)^1b_{n-k}b_{n-k-2}\frac{\theta_{j-2}^{'}}{(-1)^{j-2}}+(-1)^0b_{n-k-1}\frac{\theta_{j-1}^{'}}{(-1)^{j-1}}
\end{dmath}
\begin{dmath}
    \frac{\theta_j^{'}}{(-1)^j}=(-1)^{j-1}b_{n-k}^{j-1}\left(a_nb_{n-k-j}-a_{n-j}b_{n-k}\right)+(-1)^{j-2}b_{n-k}^{j-2}b_{n-k(j-1)}\frac{\theta_1^{'}}{(-1)^1}+(-1)^{j-3}b_{n-k}^{j-3}b_{n-k-(j-2)}\frac{\theta_2^{'}}{(-1)^2}+\dots+(-1)^2b_{n-k}^2b_{n-k-3}\frac{\theta_{j-2}^{'}}{(-1)^{j-3}}+(-1)^1b_{n-k}b_{n-k-2}\frac{\theta_{j-2}}{(-1)^{j-2}}+(-1)^0b_{n-k-1}\frac{\theta_{j-1}^{'}}{(-1)^{j-1}}
\end{dmath}
\begin{dmath}
    \frac{\theta_j^{'}}{(-1)^j}=(-1)^{j-1}b_{n-k}^{j-1}a_nb_{n-k-j}-(-1)^{j-1}b_{n-k}^{j-1}a_{n-j}b_{n-k}+(-1)^{j-2}b_{n-k}^{j-2}b_{n-k-(j-1)}\frac{\theta_1^{'}}{(-1)^1}+(-1)^{j-3}b_{n-k}^{j-3}b_{n-k-(j-2)}\frac{\theta_2^{'}}{(-1)^2}+\dots+(-1)^2b_{n-k}^2b_{n-k-3}\frac{\theta_{j-3}^{'}}{(-1)^{j-3}}+(-1)^1b_{n-k}b_{n-k-2}\frac{\theta{j-2}}{(-1)^{j-2}}+(-1)^0b_{n-k-1}\frac{\theta_{j-1}^{'}}{(-1)^{j-1}}
\end{dmath}
Since \(a_n=\theta_0^{'}\), we have
\begin{dmath}
    \frac{\theta_j^{'}}{(-1)^j}=(-1)^{j-1}b_{n-k}^{j-1}b_{n-k-j}\theta_0^{'}-(-1)^{j-1}b_{n-k}^{j-1}a_{n-j}b_{n-k}+(-1)^{j-2}b_{n-k}^{j-2}b_{n-k-(j-1)}\frac{\theta_1^{'}}{(-1)^1}+(-1)^{j-3}b_{n-k}^{j-3}b_{n-k-(j-2)}+\dots+(-1)^2b_{n-k-2}^2b_{n-k-3}\frac{\theta_{j-3}^{'}}{(-1)^{j-3}}+(-1)^1b_{n-k}b_{n-k-2}\frac{\theta_{j-2}^{'}}{(-1)^{j-2}}+(-1)^0b_{n-k-1}\frac{\theta_{j-1}^{'}}{(-1)^{j-1}}
\end{dmath}
\begin{dmath}
    \frac{\theta_j^{'}}{(-1)^j}=(-1)^jb_{n-k}^{j-1}a_{n-j}b_{n-k}+(-1)^{j-1}b_{n-k}^{j-1}b_{n-k-j}\theta_0^{'}+(-1)^{j-3}b_{n-k}^{j-2}b_{n-k-(j-1)}\theta_1^{'}+(-1)^{j-5}b_{n-k}^{j-3}b_{n-k-(j-2)}\theta_2^{'}+\dots+(-1)^{5-j}b_{n-k}^2b_{n-k-3}\theta_{j-3}^{'}+(-1)^{3-j}b_{n-k}b_{n-k-2}\theta_{j-2}^{'}+(-1)^{1-j}b_{n-k-1}\theta_{j-1}^{'}
\end{dmath}
\begin{dmath}
    \frac{\theta_j^{'}}{(-1)^j}=(-1)^jb_{n-k}^{j-1}a_{n-j}b_{n-k}+(-1)^{j+1}b_{n-k}^{j-1}b_{n-k-j}\theta_0^{'}+(-1)^{j+1}b_{n-k}^{j-2}b_{n-k-(j-1)}\theta_1^{'}+(-1)^{j+1}b_{n-k}^{j-3}b_{n-k-(j-2)}\theta_2^{'}+\dots+(-1)^{j+1}b_{n-k}^2b_{n-k-3}\theta_{j-3}^{'}+(-1)^{j+1}b_{n-k}b_{n-k-2}\theta_{j-2}^{'}+(-1)^{j+1}b_{n-k-1}\theta_{j-1}^{'}
\end{dmath}
Multiply both sides by \((-1)^j\),
\begin{dmath}
    \theta_j^{'}=(-1)^2b_{n-k}^{j-1}a_{n-j}b_{n-k}+(1)^{2j+1}b_{n-k}^{j-1}b_{n-k-j}\theta_0^{'}+(-1)^{2j+1}b_{n-k}^{j-2}b_{n-k-(j-1)}\theta_1^{'}+\dots+(-1)^{2j+1}b_{n-k}^2b_{n-k-3}\theta_{j-3}^{'}+(-1)^{2j+1}b_{n-k}b_{n-k-2}\theta_{j-2}^{'}+(-1)^{2j+1}b_{n-k-1}\theta_{j-1}^{'}
\end{dmath}
\begin{dmath}
    \theta_j^{'}=a_{n-j}b_{n-k}^j-\left(b_{n-k-j}b_{n-k}^{j-1}\theta_0^{'}+b_{n-k-(j-1)}b_{n-k}^{j-2}\theta_1^{'}+\dots+b_{n-k-2}b_{n-k}\theta_{j-2}^{'}+b_{n-k-1}\theta_{j-1}^{'}\right)
\end{dmath} 
as desired.
\( \\ \\ \)
\section{Examples}
\subsection{Oblique Asymptote}
For the linear function 
\begin{dmath}
    f(x)=\frac{{a_nx^n}+a_{n-1}x^{n-1} + \cdots + a_1x+a_0}{b_{n-1}x^{n-1}+b_{n-2}x^{n-2}+\cdots + b_1x+b_{0}}
\end{dmath} 
the formula for its oblique asymptote would be
\begin{dmath}
    y=\frac{a_n}{b_{n-1}}x-\frac{
    \begin{vmatrix}
        a_n & a_{n-1} \\
        b_{n-1} & b_{n-2}
    \end{vmatrix}}{b_{n-1}^2}
\end{dmath} where \(b_i=0\) for \(i<0\). \\
\\
\textbf{Example 4.1} Referring back to the function from Example 1.2:
\begin{dmath}
    \displaystyle f(x)=\frac{5x^3+13x^2+3x+9} {4x^2+5x+7}
\end{dmath}
The corresponding asymptote would be
\begin{dmath}
    y=\frac {5}{4}x-\frac{
    \begin{vmatrix}
        5 & 13 \\
        4 & 5
    \end{vmatrix}}{4^2}
\end{dmath}
\begin{dmath}
    \displaystyle y=\frac{5}{4}x+\frac{27}{16}
\end{dmath}
\subsection{Curvilinear Asymptote}
For the quadratic function 
\begin{dmath}
    f(x)=\frac{{a_nx^n}+a_{n-1}x^{n-1} + \cdots + a_1x+a_0}{b_{n-2}x^{n-2}+b_{n-3}x^{n-3}+\cdots + b_1x+b_{0}}
\end{dmath} 
the formula for its curvilinear asymptote becomes
\begin{dmath}
    y=\frac{a_n}{b_{n-2}}x^2-\frac{
    \begin{vmatrix}
        a_n & a_{n-1} \\
        b_{n-2} & b_{n-3}
    \end{vmatrix}}{b_{n-2}^2}x +\frac{
    \begin{vmatrix}
        a_n & a_{n-1} & a_{n-2} \\
        b_{n-2} & b_{n-3} & b_{n-4} \\
        0 & b_{n-2} & b_{n-3}
    \end{vmatrix}}{b_{n-2}^3}
\end{dmath}
where \(b_i=0\) for \(i<0\). \\
\\
\textbf{Example 4.2} The curvilinear asymptote of the function
\begin{dmath}
\displaystyle f(x)=\frac {x^4-2x^3+3x-9}{2x^2-5}\\ 
\end{dmath}
would be
\begin{dmath}
    y=\frac{1}{2}x^2-\frac{
    \begin{vmatrix}
        1 & -2 \\
        2 & 0
    \end{vmatrix}}{2^2}x +\frac{
    \begin{vmatrix}
        1 & -2 & 0 \\
        2 & 0 & -5 \\
        0 & 2 & 0
    \end{vmatrix}}{2^3} 
\end{dmath}
\begin{dmath}
    \displaystyle y=\frac {x^2}{2} - x+ \frac {5}{4}
\end{dmath}
\( \\ \)
\textbf{Example 4.3} Referring back to the function from Example 1.1:
\begin{dmath}
\displaystyle f(x)=\frac {8x^3+7}{x-4} \\ 
\end{dmath}
the corresponding quadratic asymptote would be
\begin{dmath}
    y=\frac{8}{1}x^2-\frac{
    \begin{vmatrix}
        8 & 0 \\
        1 & -4
    \end{vmatrix}}{1^2}x+\frac{
    \begin{vmatrix}
        8 & 0 & 0 \\
        1 & -4 & 0 \\
        0 & 1 & -4
    \end{vmatrix}}{1^3} 
\end{dmath}
\begin{dmath}
    \displaystyle y= 8x^2+32x+128 
\end{dmath}
\( \\ \\ \)
\section {Acknowledgements}
The authors would like to thank two Greek students, Aris Kapras and Nikos Parlitsis, who inspired the development of the formula presented in this paper. In addition, many thanks to the people on Quora for the useful discussions, and special thanks to Brian Rago for his insight into the general form of the formula.
\( \\ \\ \)
\section{Bibliography}
\([1]\) Britannica, T., Editors of Encyclopaedia (Invalid Date). Asymptote. Encyclopedia Britannica. https://www.britannica.com/science/asymptote  \\
\\
\([2]\) Sendra, J.R., Winkler, F., Perez-Diaz, S. (2007). Rational Algebraic Curves: A Computer Algebra Approach. Series: Algorithms and Computation in Mathematics. Vol. 22. Springer Verlag \\
\\
\([3]\) Stewart, J., Clegg, D., Watson, S., Cengage Learning. (2021). Calculus. Australia [i pozostałe: Cengage Learning.\\
\\
\([4]\) Thomas, G. B., Finney, R. L., Weir, M. D. (1996). Calculus and Analytic Geometry. Reading, Mass: Addison-Wesley Pub. Co. \\

\end{document}